\documentclass[10pt,letterpaper]{article}

\makeatletter

\usepackage{latexsym}

\makeatother

\usepackage{blindtext,titlefoot}
\usepackage{authblk}
\usepackage{marvosym}
\usepackage[utf8]{inputenc}
\usepackage[english]{babel}
\usepackage{amsmath}
\usepackage{amsfonts}
\usepackage{amssymb}
\usepackage{makeidx}
\usepackage{graphicx}
\usepackage{lmodern}
\usepackage{kpfonts}
\usepackage{dsfont}
\usepackage{cancel}
\usepackage{amsthm} 
\usepackage[left=2cm,right=2cm,top=2cm,bottom=2cm]{geometry}
\usepackage{subcaption}
\usepackage{multirow}
\usepackage{float}
\usepackage{url}
\usepackage{breakurl}
\usepackage[breaklinks]{hyperref}
\usepackage{multicol}
\usepackage{mathtools, cuted}
\usepackage{lipsum}
\usepackage{booktabs}
\usepackage{comment}
\usepackage{cite}

\usepackage{fancyhdr}
\fancypagestyle{firstpage}{\fancyhf{}
\setcounter{page}{1}
\rfoot{\thepage}
}

\providecommand{\nset}[1]{
\mathbb{#1}
}
\providecommand{\set}[1]{
\left\{#1\right\}
}
\providecommand{\com}[1]{``#1"}

\providecommand{\ifr}[5]{
{}^{#1}_{#2}{#3}_{#4}^{#5}
}
\providecommand{\gam}[1]{
\Gamma\left(#1 \right)
}

\providecommand{\norm}[1]{
\left\lVert #1 \right\rVert
}
\providecommand{\abs}[1]{
\left\lvert #1 \right\rvert
}
\providecommand{\ds}[1]{
\displaystyle #1
}
\providecommand{\der}[3]{
\dfrac{#1^{#3} }{ #1 #2^{#3}}
}

\providecommand{\re}[1]{
\hbox{Re}\left(#1 \right)
}
\providecommand{\im}[1]{
\hbox{Im}\left(#1 \right)
}
\providecommand{\rnd}[2]{
\hbox{Rnd}_#2\left(#1\right)
}

\newtheorem{theorem}{ Theorem}[section]
\newtheorem{definition}[theorem]{Definition}
\newtheorem{proposition}[theorem]{Proposition}
\newtheorem{corollary}[theorem]{Corollary}
\newtheorem{example}[theorem]{Example}

\usepackage{enumitem} 
\setlist[itemize]{noitemsep} 

\usepackage{titlesec} 
\titleformat{\section}[block]{\large\bfseries\scshape\centering}{\thesection.}{1em}{} 
\titleformat{\subsection}[block]{\large\bfseries\scshape\centering}{\thesubsection.}{1em}{}
\titleformat{\subsubsection}[block]{\large\bfseries\scshape\centering}{\thesubsubsection.}{1em}{} 

\title{\Huge\bfseries Fractional Newton-Raphson Method \\ Accelerated with Aitken's Method}

\author[,a,$\star$]{\normalsize A. Torres-Hernandez  \footnote{E-mail: anthony.torres@ciencias.unam.mx;  ORCID: 0000-0001-6496-9505}}
\affil[a]{\normalsize Department of Physics, Faculty of Science - UNAM, Mexico}

\author[,b,$\star$]{\normalsize F. Brambila-Paz \footnote{E-mail: fernando.brambila@ciencias.unam.mx; ORCID: 0000-0001-7896-6460} }
\affil[b]{\normalsize Department of Mathematics, Faculty of Science - UNAM, Mexico}

\author[,b]{\normalsize U. Iturrarán-Viveros \footnote{E-mail: ursula@ciencias.unam.mx; ORCID: 0000-0001-6410-7542}}

\author[,a]{\normalsize R. Caballero-Cruz  \footnote{E-mail: rcaballero@ciencias.unam.mx; ORCID: 0000-0003-4470-9399}}

\date{}

\begin{document}

\maketitle

\thispagestyle{firstpage}

\maketitle\unmarkedfntext{$\star$ Corresponding authors}

\begin{abstract}

In the following document, we present a way to obtain the order of convergence of the Fractional Newton-Raphson (F N-R) method, which seems to have an order of convergence at least linearly for the case in which the order $\alpha$ of the derivative is different from one. A simplified way of constructing the Riemann-Liouville (R-L) fractional operators, fractional integral and fractional derivative, is presented along with examples of its application on different functions. Furthermore, an introduction to the Aitken's method is made and it is explained why it has the ability to accelerate the convergence of the iterative methods, to finally present the results that were obtained when implementing the Aitken's method in the F N-R method.

\textbf{Keywords:} Newton-Raphson Method, Fractional Calculus, Fractional Derivative, Aitken’s Method.
\end{abstract}

\section{Fixed Point Method}

A classic problem in mathematics, which is of common interest in physics and engineering, is finding the set of zeros of a function $f:\Omega \subset \nset{R}^n \to \nset{R}^n$, that is,

\begin{eqnarray}\label{eq:1-001}
\begin{array}{c}
\set{\xi \in \Omega \ : \ \norm{f(\xi)}=0},
\end{array}
\end{eqnarray}

where $\norm{ \ \cdot \ }: \nset{R}^n \to \nset{R}$ denotes any vector norm. Although finding the zeros of a function may seem like a simple problem, in general it involves solving an \textbf{algebraic equation system} as follows

\begin{eqnarray}\label{eq:1-002}
\left\{
\begin{array}{c}
\left[f\right]_1(x)=0\\
\left[f\right]_2(x)=0\\
\vdots \\
\left[f\right]_n(x)=0
\end{array}\right.,
\end{eqnarray}

where $[f]_k: \nset{R}^n \to \nset{R}$ denotes the $k$-th component of the function $f$. It should be noted that the system of equations \eqref{eq:1-002} may represent a \textbf{linear system} or a \textbf{nonlinear system}, and in general, it is necessary to use numerical methods of the iterative type to solve it. Let  $\Phi:\nset{R}^n \to \nset{R}^n$ be a function, it is possible to build a sequence $\set{x_i}_{i=0}^\infty$  by defining the following iterative method

\begin{eqnarray}\label{eq:1-003}
x_{i+1}:=\Phi(x_i),
\end{eqnarray}

if it is fulfills that $x_i\to \xi\in \nset{R}^n$, and if the function $\Phi$ is continuous around $\xi$, we obtain that

\begin{eqnarray}\label{eq:1-004}
\xi=\lim_{i\to \infty}x_{i+1}=\lim_{i\to \infty}\Phi(x_i)=\Phi\left(\lim_{i\to \infty}x_i \right)=\Phi(\xi),
\end{eqnarray}

the above result is the reason by which the method \eqref{eq:1-003} is known as the \textbf{fixed point method}. Furthermore, the function $\Phi$ is called an \textbf{iteration function}. To understand the nature of the convergence of the iteration function $\Phi$, the following definition is necessary \cite{plato2003concise}:

\begin{definition}
Let $\Phi:\nset{R}^n \to \nset{R}^n$  be an iteration function. The method \eqref{eq:1-003} for determining $\xi\in \nset{R}^n$ is called \textbf{(locally) convergent}, if there exists $\delta>0$ such that for every initial value

\begin{eqnarray*}
x_0\in B(\xi;\delta):=\set{y\in \nset{R}^n \ : \  \norm{y-\xi}<\delta},
\end{eqnarray*}

it is fulfills that

\begin{eqnarray}\label{eq:1-005}
\lim_{i \to \infty}\norm{x_i-\xi}\to 0 & \Rightarrow & \lim_{i\to \infty}x_i=\xi.
\end{eqnarray}

\end{definition}

If we have a function $f:\Omega \subset \nset{R}^n \to \nset{R}^n$, for which we want to determine the set \eqref{eq:1-001}, in general it is possible to write an iteration function 
$\Phi$ as follows \cite{burden2002analisis}

\begin{eqnarray*}
\Phi(x)=x-A(x)f(x),
\end{eqnarray*}

where $A(x)$  is a matrix, which is given as follows

\begin{eqnarray*}
A(x):=\left([A]_{jk}(x) \right)=\begin{pmatrix}
[A]_{11}(x)&[A]_{12}(x)& \cdots &[A]_{1n}(x)\\
[A]_{21}(x)&[A]_{22}(x)&\cdots& [A]_{2n}(x)\\
\vdots & \vdots & \ddots & \vdots \\
[A]_{n1}(x)&[A]_{n2}(x)&\cdots &[A]_{nn}(x)
\end{pmatrix},
\end{eqnarray*}

with $[A]_{jk}(x):\nset{R}^n \to \nset{R} \ \forall j,k\leq n$. It is necessary to mention that the matrix  $ A (x) $ is determined according to the order of convergence desired.

\subsection{Order of Convergence}

Before continuing, it is necessary to define the order of convergence of an iteration function $\Phi$ \cite{plato2003concise}:

\begin{definition}
Let $ \Phi: \Omega \subset \nset{R}^ n \to \nset{R}^ n $ be an iteration function with a fixed point $ \xi \in \Omega $. Then the method \eqref{eq:1-003} is called  \textbf{(locally) convergent of (at least) order $ \boldsymbol{p} $} ($ p \geq 1 $), if there are exists $ \delta> 0 $  and $ C $ a non-negative constant,  with $ C <1 $ if $ p = 1 $, such that for any initial value $ x_0 \in B (\xi; \delta) $ it is fulfills that

\begin{eqnarray}\label{eq:1-006}
\norm{x_{k+1}-\xi}\leq C \norm{x_k-\xi}^p, & k=0,1,2,\cdots,
\end{eqnarray}

where $ C $ is called convergence factor.

\end{definition}

The order of convergence is usually related to the speed at which the sequence generated by \eqref{eq:1-003} converges. For the particular case $ p = 1 $ it is said that the method \eqref{eq:1-003} has an \textbf{order of convergence (at least) linear}, and for the case $p=2$  it is said that the method \eqref{eq:1-003} has an \textbf{order of convergence (at least) quadratic}. The following theorem, allows characterizing the order of convergence of an iteration function $ \Phi $ with its derivatives \cite{plato2003concise,stoer2013,burden2002analisis,torres2020reduction}. Before continuing,  we need to consider the following multi-index notation. Let $\nset{N}_0$ be the set $\nset{N}\cup\set{0}$, if $\gamma \in \nset{N}_0^n$, then

\begin{eqnarray}
\left\{
\begin{array}{l}
 \gamma!:= \ds\prod_{k=1}^n [\gamma]_k ! \vspace{0.1cm}\\
 \abs{\gamma}:= \ds \sum_{k=1}^n [\gamma]_k\vspace{0.1cm}\\
 x^\gamma:= \ds \prod_{k=1}^n [x]_k^{[\gamma]_k}\vspace{0.1cm}\\
\der{\partial}{x}{\gamma}:= \dfrac{\partial^{\abs{\gamma}}}{\partial [x]_1^{[\gamma]_1}\partial [x]_2^{[\gamma]_2}\cdots \partial [x]_n^{[\gamma]_n} }
\end{array}\right. .
\end{eqnarray}

\begin{theorem}\label{teo:1-001}
Let $ \Phi: \Omega \subset \nset{R}^n \to \nset {R}^n $ be an iteration function with a fixed point $ \xi \in \Omega $. Assuming that $\Phi $ is $ p$-times differentiable in $ \xi $ for some $ p \in \nset{N} $, and furthermore

\begin{eqnarray}\label{eq:1-007}
\left\{
\begin{array}{cc}
\ds  \dfrac{\partial^\gamma [\Phi]_k(\xi) }{ \partial x^\gamma}=0, \ \forall k\geq 1 \mbox{ and } \forall  \abs{\gamma}<p, & \mbox{if }p\geq 2 \vspace{0.1cm}\\
\ds \norm{\Phi^{(1)}(\xi)}<1, & \mbox{if }p=1
\end{array}\right.,
\end{eqnarray}

where $\Phi^{(1)}$ denotes the \textbf{Jacobian matrix} of the function $\Phi$, then $ \Phi $ is (locally) convergent of (at least) order $ p $.

\begin{proof}

Let $\Phi:\nset{R}^n \to \nset{R}^n$ be an iteration function, and let $\set{\hat{e}_k}_{k=1}^n$ be the canonical basis of $\nset{R}^n$. Considering the following index notation (Einstein notation)

\begin{eqnarray*}
\Phi(x)=\sum_{k=1}^n [\Phi]_k(x)\hat{e}_k: = [\Phi]_k(x)\hat{e}_k=\hat{e}_k[\Phi]_k(x),
\end{eqnarray*}

and using the Taylor series expansion of a vector-valued function in multi-index notation, we obtain two cases:

\begin{itemize}

\item[i)]  Case $p\geq 2:$

\begin{align*}
\Phi(x_i)
=& \ds  \Phi(\xi)+  \sum_{\abs{\gamma} =1}^p \dfrac{1}{\gamma !}\hat{e}_k\dfrac{\partial^\gamma [\Phi]_k(\xi) }{ \partial x^\gamma}   (x_i-\xi)^\gamma   + \hat{e}_k[o]_k\left(\max_{\abs{\gamma}=p}  \set{ (x_i-\xi)^\gamma }\right) \\
=& \ds \Phi(\xi)+ \sum_{m =1}^p   \left( \sum_{\abs{\gamma} =m}\dfrac{1}{\gamma !} \hat{e}_k \dfrac{\partial^\gamma [\Phi]_k(\xi) }{ \partial x^\gamma}   (x_i-\xi)^\gamma \right)   + \hat{e}_k[o]_k\left(\max_{\abs{\gamma}=p}  \set{ (x_i-\xi)^\gamma }\right),
\end{align*}

then

\begin{align*}
\norm{\Phi(x_i)-\Phi(\xi)}&\leq  \ds \sum_{m =1}^p   \left( \sum_{\abs{\gamma} =m}\dfrac{1}{\gamma !}\norm{ \hat{e}_k \dfrac{\partial^\gamma [\Phi]_k(\xi) }{ \partial x^\gamma}   (x_i-\xi)^\gamma }\right)  +  \norm{\hat{e}_k  [o]_k\left(\max_{\abs{\gamma}=p}  \set{ (x_i-\xi)^\gamma }\right)} \\
&\leq   \ds \sum_{m =1}^p   \left( \sum_{\abs{\gamma} =m}\dfrac{1}{\gamma !}\norm{ \dfrac{\partial^\gamma [\Phi]_k(\xi) }{ \partial x^\gamma}  \hat{e}_k   }\right)\norm{x_i-\xi}^m+  o\left( \norm{x_i-\xi}^p \right), 
\end{align*}

assuming that $\xi$ is a fixed point of $\Phi$ and that $\dfrac{\partial^\gamma [\Phi]_k(\xi) }{ \partial x^\gamma}=0 \ \forall k\geq 1$ and $  \forall \abs{\gamma}<p$ is fulfilled, the previous expression implies that

\begin{eqnarray*}
\dfrac{\norm{\Phi(x_i)-\Phi(\xi)}}{\norm{x_i-\xi}^p}=\dfrac{\norm{x_{i+1}-\xi}}{\norm{x_i-\xi}^p}\leq\sum_{\abs{\gamma} =p}\dfrac{1}{\gamma !}\norm{ \dfrac{\partial^\gamma [\Phi]_k(\xi) }{ \partial x^\gamma}  \hat{e}_k   } +\dfrac{o\left(\norm{x_i-\xi}^p \right)}{\norm{x_i-\xi}^p},
\end{eqnarray*}

therefore

\begin{eqnarray*}
\lim_{i\to \infty} \dfrac{\norm{x_{i+1}-\xi}}{\norm{x_i-\xi}^p}\leq  \sum_{\abs{\gamma} =p}\dfrac{1}{\gamma !}\norm{ \dfrac{\partial^\gamma [\Phi]_k(\xi) }{ \partial x^\gamma}  \hat{e}_k   },
\end{eqnarray*}

as a consequence, if the sequence $\set{x_i}_{i=0}^\infty$ generated by \eqref{eq:1-003} converges to $\xi$, there exists a value $k>0$ such that

\begin{eqnarray*}
\norm{x_{i+1}-\xi}\leq \left( \sum_{\abs{\gamma} =p}\dfrac{1}{\gamma !}\norm{ \dfrac{\partial^\gamma [\Phi]_k(\xi) }{ \partial x^\gamma}  \hat{e}_k   }\right)\norm{x_i-\xi}^p, & \forall i\geq k,
\end{eqnarray*}

then $ \Phi $ is (locally) convergent of (at least) order $ p $.

 \item[ii)] Case $p=1:$

\begin{align*}
\Phi(x_i)=& \ds \Phi(\xi)+  \sum_{\abs{\gamma} =1}\dfrac{1}{\gamma !} \hat{e}_k \dfrac{\partial^\gamma [\Phi]_k(\xi) }{ \partial x^\gamma}   (x_i-\xi)^\gamma    + \hat{e}_k[o]_k\left(\max_{\abs{\gamma}=1}  \set{ (x_i-\xi)^\gamma }\right)\\
=&\Phi(\xi)+\Phi^{(1)}(x_i)(x_i-\xi)+\hat{e}_k[o]_k\left(\max_{\abs{\gamma}=1}  \set{ (x_i-\xi)^\gamma }\right),
\end{align*}

then

\begin{align*}
\norm{\Phi(x_i)-\Phi(\xi)}\leq  \norm{\Phi^{(1)}(\xi)} \norm{x_i-\xi }+   o\left( \norm{x_i-\xi} \right) ,
\end{align*}

assuming that $\xi$ is a fixed point of $\Phi$, the previous expression implies that

\begin{eqnarray*}
\dfrac{\norm{\Phi(x_i)-\Phi(\xi)}}{\norm{x_i-\xi}}=\dfrac{\norm{x_{i+1}-\xi}}{\norm{x_i-\xi}}\leq \norm{\Phi^{(1)}(\xi)} +\dfrac{o\left(\norm{x_i-\xi} \right)}{\norm{x_i-\xi}},
\end{eqnarray*}

therefore

\begin{eqnarray*}
\lim_{i\to \infty} \dfrac{\norm{x_{i+1}-\xi}}{\norm{x_i-\xi}}\leq\norm{\Phi^{(1)}(\xi)},
\end{eqnarray*}

as a consequence, if the sequence $\set{x_i}_{i=0}^\infty$ generated by \eqref{eq:1-003} converges to $\xi$, there exists a value $k>0$ such that

\begin{eqnarray*}
\norm{x_{i+1}-\xi}\leq \norm{\Phi^{(1)}(\xi)}\norm{x_i-\xi}, & \forall i\geq k,
\end{eqnarray*}

considering $m\geq 1$, from the previous inequality we obtain that

\begin{align*}
\norm{x_{i+m}-\xi}\leq &\norm{\Phi^{(1)}(\xi)}\norm{x_{i+m-1}-\xi}\leq \norm{\Phi^{(1)}(\xi)}^2\norm{x_{i+m-2}-\xi} \leq \cdots \leq\norm{\Phi^{(1)}(\xi)}^{m}\norm{x_{i}-\xi},
\end{align*}

and assuming that $\norm{\Phi^{(1)}(\xi)}<1$ is fulfilled

\begin{eqnarray*}
\lim_{m \to \infty}\norm{x_{i+m}-\xi}\leq \lim_{m \to \infty}\norm{\Phi^{(1)}(\xi)}^{m}\norm{x_{i}-\xi} \to 0,
\end{eqnarray*}

then $ \Phi $ is (locally) convergent of order (at least) linear.

\end{itemize}

\end{proof}
\end{theorem}

The following corollary follows from the previous theorem

\begin{corollary}\label{cor:1-001}
Let $\Phi:\nset{R}^n \to \nset{R}^n$ be an iteration function. If $\Phi$ defines a sequence $\set{x_i}_{i=0}^\infty$ such that $x_i\to \xi$, and if the following condition is true

\begin{eqnarray}\label{eq:1-008}
\lim_{x\to \xi}\norm{\Phi^{(1)}(x)}\neq 0,
\end{eqnarray}

then $\Phi$ has an order of convergence (at least) linear  in $B(\xi;\delta)$.
\end{corollary}

\section{Newton-Raphson Method}

We begin this section by considering the following proposition   \cite{torres2020reduction,torreshern2020} :

\begin{proposition}\label{prop:2-001}
Let $f:\Omega \subset \nset{R}^n\to \nset{R}^n$ be a function with a value $\xi\in \Omega$ such that $\norm{f(\xi)}=0$, and let $\Phi:\nset{R}^n \to \nset{R}^n$ be an iteration function as follows

\begin{eqnarray}\label{eq:2-001}
\Phi(x)=x-A(x)f(x),
\end{eqnarray}

with $A(x)$ a matrix. If the following condition it is fulfills

\begin{eqnarray}\label{eq:2-002}
\lim_{x\to \xi}A(x)=\left(f^{(1)}(\xi) \right)^{-1},
\end{eqnarray}

where $f^{(1)}$ denotes the Jacobian matrix of the function $f$,  which is defined as follows \cite{ortega1990numerical}

\begin{eqnarray}\label{eq:2-003}
f^{(1)}(x):=\left([f]^{(1)}_{jk}(x) \right) =\left( \partial_k[f]_j(x) \right),
\end{eqnarray}

where

\begin{eqnarray*}
[f]^{(1)}_{jk}(x)=\partial_k[f]_j(x):= \dfrac{\partial }{\partial [x]_k}[f]_j(x), &1\leq j,k\leq n,
\end{eqnarray*}

then the iteration function $\Phi$, fulfills a necessary (but not sufficient) condition to be (locally) convergent of order (at least) quadratic in $B(\xi;\delta)$.

\begin{proof}
From the \textbf{Theorem \ref{teo:1-001}}, we have that an iteration function has an order of convergence (at least) quadratic if it fulfills the following condition

\begin{eqnarray*}
\lim_{x\to \xi}\dfrac{\partial [\Phi]_k(x)}{\partial [x]_j}=0, & \forall j,k\leq n,
\end{eqnarray*}

which may be written equivalently as follows

\begin{eqnarray}\label{eq:2-004}
 \lim_{x\to \xi}\norm{\Phi^{(1)}(x)}=0.
\end{eqnarray}

Then, we may assume that we have a function $f(x):\Omega \subset \nset{R}^n \to \nset{R}^n$  with a zero $\xi \in \Omega$, such that all of its first partial derivatives are defined in $ \xi $, and taking the iteration function $ \Phi $ given by \eqref{eq:2-001}, the $ k $-th component of the iteration function may be written as

\begin{eqnarray*}
[\Phi]_k(x)=[x]_k-\sum_{j=1}^n[A]_{kj}(x)[f]_j(x),
\end{eqnarray*}

then

\begin{align*}
\partial_l [\Phi]_k(x)=&\delta_{lk}-\sum_{j=1}^n \big{(} [A]_{kj}(x)\partial_l [f]_j(x)+\left(\partial_l [A]_{kj}(x) \right)[f]_j(x) \big{)} \\
=&\delta_{kl}-\sum_{j=1}^n \left( [A]_{kj}(x)[f]^{(1)}_{jl}(x)+\left(\partial_l [A]_{kj}(x) \right)[f]_j(x) \right),
\end{align*}

where $ \delta_{kl} $ is the Kronecker delta, which is defined as

\begin{eqnarray*}
\delta_{kl}=\delta_{lk}=\left\{
\begin{array}{cc}
1,& \mbox{si }l=k\\
0,& \mbox{si }l\neq k
\end{array}\right..
\end{eqnarray*}

Assuming that condition \eqref{eq:2-004} it is fulfills

\begin{eqnarray*}
\partial_l [\Phi]_k(\xi)=\delta_{kl}-\sum_{j=1}^n [A]_{kj}(\xi)[f]^{(1)}_{jl}(\xi)=0 & \Rightarrow & \sum_{j=1}^n [A]_{kj}(\xi)[f]^{(1)}_{jl}(\xi)=\delta_{kl}, \hspace{0.1cm} \forall l,k\leq n,
\end{eqnarray*}

then the above expression may be written in matrix form as follows

\begin{eqnarray*}
A(\xi) f^{(1)}(\xi)=I_n & \Rightarrow & A(\xi)= \left(f^{(1)}(\xi)\right)^{-1},
\end{eqnarray*}

where  $ I_n $ denotes the identity matrix of $ n \times n $. Then any matrix $ A (x) $ that fulfills the following condition

\begin{eqnarray*}
\lim_{x\to \xi}A(x)= \left(f^{(1)}(\xi)\right)^{-1},
\end{eqnarray*}

guarantees that exists $ \delta> 0 $, such that iteration function $ \Phi $ given by \eqref{eq:2-001}, fulfills a necessary (but not sufficient) condition to be (locally) convergent of  order (at least) quadratic in $ B (\xi; \delta) $.

\end{proof}
\end{proposition}

The following fixed point method may be obtained from the previous proposition

\begin{eqnarray}\label{eq:2-005}
x_{i+1}:=\Phi(x_i)=x_i-\left(f^{(1)}(x_i) \right)^{-1}f(x_i) , & i=0,1,2,\cdots,
\end{eqnarray}

which is known as \textbf{Newton-Raphson method}, also known as Newton's method \cite{ortega1970iterative}. Given the condition \eqref{eq:2-002}, it could be wrongly considered that the Newton-Raphson method always has an order of convergence (at least) quadratic, but as mentioned in the \textbf{Proposition \ref{prop:2-001}}, the form of the iteration function \eqref{eq:2-005} is not sufficient to guarantee this order of convergence. This occurs because even if the condition \eqref{eq:2-002} it is fulfills, the order of convergence  becomes conditioned by the way in which the function $f$ is constituted, for example for the one variable case, if the function $f$ has a root $\xi$, with a certain algebraic multiplicity $ m\geq 2$, that is, 

\begin{eqnarray*}
\begin{array}{cc}
f(x)=(x-\xi)^mg(x), & g(\xi)\neq 0,
\end{array}
\end{eqnarray*}

the Newton-Raphson method presents an order of convergence at least linear \cite{plato2003concise}, the aforementioned may be seen in the following proposition:

 \begin{proposition}\label{prop:2-002}
Let $f:\Omega \subset \nset{R} \to \nset{R}$ be a function with a zero $\xi \in \Omega$. Then the iteration function $\Phi$ of the Newton-Raphson method, given by \eqref{eq:2-005}, fulfills the following condition:

\begin{eqnarray}\label{eq:2-006}
\abs{x_{i+1}- \xi}\leq   \dfrac{\abs{\Phi^{(p)}(\xi)}}{p !}\abs{x_i-\xi}^p,
\end{eqnarray}

where

\begin{eqnarray}
p=\left\{
\begin{array}{cl}
1 ,& \mbox{ if } f(x)=(x-\xi)^mg(x)  \vspace{0.1cm}\\
2, & \mbox{ if } f^{(1)}(\xi)\neq 0, \mbox{ and }f(x)\neq (x-\xi)^mg(x) \vspace{0.1cm}\\
3,& \mbox{ if } f^{(1)}(\xi)\neq 0, \hspace{0.1cm} f^{(2)}(\xi)=0, \mbox{ and }f(x)\neq (x-\xi)^mg(x) \vspace{0.1cm} \\
4,&\mbox{ if } f^{(1)}(\xi)\neq 0, \hspace{0.1cm} f^{(2)}(\xi)=0,\hspace{0.1cm} f^{(3)}(\xi)=0 \mbox{ and }f(x)\neq (x-\xi)^mg(x) 
\end{array}\right. ,
\end{eqnarray}

with  $g(\xi)\neq 0$ and $m\geq 2$.

\begin{proof}
Considering that the form of the function $f$ is not explicitly determined, it is possible to consider two possibilities:

\begin{itemize}
\item[i)] Assuming the function may be written as $f(x)=(x-\xi)^mg(x)$  with  $g(\xi)\neq 0$  and $m\geq 2$, then

\begin{eqnarray*}
f^{(1)}(x)=(x-\xi)^{m-1}\left[(x-\xi) g^{(1)}(x)+mg(x)\right],
\end{eqnarray*}

as a consequence, the iteration function of N-R method takes the following form

\begin{eqnarray*}
\Phi(x)= x-(x-\xi)h(x)g(x),
\end{eqnarray*}

with

\begin{eqnarray*}
h(x)=\left[(x-\xi) g^{(1)}(x)+mg(x)\right]^{-1},
\end{eqnarray*}

then

\begin{align*}
\Phi^{(1)}(x)= 1- h(x)\left[(x-\xi)g^{(1)}(x)+g(x)  \right]-(x-\xi)h^{(1)}(x)g(x),
\end{align*}

where

\begin{eqnarray*}
h^{(1)}(x)=-\left[(x-\xi) g^{(1)}(x)+mg(x)\right]^{-2}\left[\left( 1+m \right) g^{(1)}(x)+(x-\xi) g^{(2)}(x) \right],
\end{eqnarray*}

therefore

\begin{align}
\lim_{x\to \xi}\abs{\Phi^{(1)}(x)}= \abs{1-h(\xi)g(\xi)}=\abs{1-\dfrac{1}{m}}<1,
\end{align}

and from the \textbf{Theorem \ref{teo:1-001}}, the Newton-Raphson method has an order of convergence at least linear, that is, fulfills the equation \eqref{eq:2-006} with $p=1$.

\item[ii)] Assuming that $f(x)\neq (x-\xi)^mg(x)$  with  $g(\xi)\neq 0$  and $m\geq 2$, the first derivative of the iteration function of Newton-Raphson method takes the following form

\begin{align*}
\Phi^{(1)}(x)=f(x)\left[\left(f^{(1)}(x) \right)^{-2}f^{(2)}(x)\right],
\end{align*}

and if it is fulfills that $f^{(1)}(\xi)\neq 0$, then

\begin{align}
\lim_{x\to \xi}\abs{\Phi^{(1)}(x)}= 0,
\end{align}

and from the \textbf{Theorem \ref{teo:1-001}}, the Newton-Raphson method has an order of convergence at least quadratic, that is, fulfills the equation \eqref{eq:2-006} with $p=2$. On other hand, the second derivative of the iteration function of Newton-Raphson method takes the following form

\begin{align*}
\Phi^{(2)}(x)= \left(f^{(1)}(x) \right)^{-1}f^{(2)}(x)+ f(x)\left[\left(f^{(1)}(x) \right)^{-2} f^{(3)}(x)-2 \left(f^{(1)}(x) \right)^{-3}\left( f^{(2)}(x) \right)^2 \right],
\end{align*}

and if it is fulfills that $f^{(1)}(\xi)\neq 0$ and $f^{(2)}(\xi)=0$, then

\begin{align}
\lim_{x\to \xi}\abs{\Phi^{(1)}(x)}=\lim_{x\to \xi}\abs{\Phi^{(2)}(x)}= 0,
\end{align}

and from the \textbf{Theorem \ref{teo:1-001}}, the Newton-Raphson method has an order of convergence at least cubic, that is, fulfills the equation \eqref{eq:2-006} with $p=3$. Finally, the third derivative of the iteration function of the Newton-Raphson method takes the following form

\begin{align*}
\Phi^{(3)}(x)=&f(x)\left[\left(f^{(1)}(x)\right)^{-2} f^{(4)}(x)\right] + 2\left(f^{(1)}(x)\right)^{-1} f^{(3)}(x)- 3 \left(f^{(1)}(x)\right)^{-2}\left(f^{(2)}(x)\right)^2 \\
& + 6 f(x) \left[\left(f^{(1)}(x)\right)^{-4}\left( f^{(2)}(x)\right)^3 - \left(f^{(1)}(x)\right)^{-3} f^{(2)}(x) f^{(3)}(x) \right] ,
\end{align*}

and if it is fulfills that $f^{(1)}(\xi)\neq 0$, $f^{(2)}(\xi)=0$ and $f^{(3)}(\xi)=0$, then

\begin{align}
\lim_{x\to \xi}\abs{\Phi^{(1)}(x)}=\lim_{x\to \xi}\abs{\Phi^{(2)}(x)}=\lim_{x\to \xi}\abs{\Phi^{(3)}(x)}= 0,
\end{align}

and from the \textbf{Theorem \ref{teo:1-001}}, the Newton-Raphson method has an order of convergence at least tetrahedral, that is, fulfills the equation \eqref{eq:2-006} with $p=4$.

\end{itemize}
\end{proof}
\end{proposition}

The previous proposition, illustrates two important points that are worth mentioning when using the N-R method to find the zeros of a function $f$:

\begin{itemize}
\item[i)] When it is not evident, unless it is explicitly specified that the function $f$ has no roots of algebraic multiplicity $m\geq 2$, technically there exists the possibility that the N-R method has an order of convergence at least linear, that is, the N-R method may fulfill the equation \eqref{eq:2-006} with $p\geq 1$.

\item[ii)] Due that the N-R method is a local iterative method, even if it proves that for a root $\xi\in \Omega$ the method has an order of convergence at least linear, this does not rule out that for the same function  $f$ it may present a higher order of convergence over the same region $\Omega$. As an example of the above, we may consider the following function

\begin{eqnarray*}
f(x)=(x-\eta)(x-\xi)^mg(x), & g(\eta)\neq g(\xi)\neq 0,
\end{eqnarray*}

with $ \eta,\xi \in \Omega$, $\abs{\eta-\xi}<\epsilon$, and  $m\geq 2$.

\end{itemize} 
 
The previous points are important, because when the N-R method is implemented in a function $f$, the zeros of the function are assumed to be unknown, and their algebraic multiplicities $m\geq 2$, in case they exist, are also unknown. With the above in mind, the following corollary is obtained, which is derived from  \textbf{Proposition \ref{prop:2-001}}, \textbf{Proposition \ref{prop:2-002}} and  \textbf{Corollary \ref{cor:1-001}}

\begin{corollary}
Let $f:\Omega \subset \nset{R}^n \to \nset{R}^n$ be a function with a zero $\xi \in \Omega$. If there exists at least a value $k>0$, and a function $g_k: \nset{R}^n \to \nset{R}$, such that,

\begin{eqnarray*}
[f]_{j_k}(x)=[(x-\xi)]_k^m g_k(x), & g_k(\xi)\neq 0,
\end{eqnarray*}

for some value $j_k$, with 

\begin{eqnarray*}
1\leq j_k,k\leq n & \mbox{ and }& m\geq 2,
\end{eqnarray*}

then the Jacobian matrix of the iteration function $\Phi$ of the N-R method, given by \eqref{eq:2-005}, fulfills that all entries in its $k$-th column are nonzero at the value $\xi$, that is,

\begin{eqnarray*}
[\Phi]_{jk}^{(1)}(\xi)\neq 0, & \forall j>0,
\end{eqnarray*}

as consequence, the N-R method has an order of convergence (at least) linear.

\end{corollary}

\section{Fractional Calculus}

The fractional calculus is a mathematical analysis branch whose applications have been increasing since the end of the XX century and beginnings of the XXI century \cite{brambila2017fractal,martinez2017applications1,martinez2017applications2}, the fractional calculus arises around 1695 due to Leibniz's notation for the derivatives of integer order

\begin{eqnarray*}
f^{(n)}(x):=\dfrac{d^n}{dx^n}f(x), & n\in \nset{N},
\end{eqnarray*}

thanks to this notation L'Hopital could ask in a letter to Leibniz about the interpretation of taking $ n = 1/2 $ in a derivative, since at that moment Leibniz could not give a physical or geometrical interpretation to this question, he simply answered L'Hopital in a letter, \com{$\dots$ is an apparent paradox of which, one day, useful consequences will be drawn} \cite{miller93}. The name of fractional calculus comes from a historical question since in this branch of mathematical analysis it is studied the derivatives and integrals of a certain order $\alpha$, with $\alpha \in \nset{R}$ or $\nset{C}$. 

Currently, the fractional calculus does not have a unified definition of what is considered a fractional derivative, because one of the conditions required to consider an expression as a fractional derivative is to recover the results of conventional calculus when the order $\alpha \to n$, with $n \in \nset{N}$ \cite{oldham74}, among the most common definitions of fractional derivatives are the  Riemann-Liouville (R-L) fractional derivative and the Caputo fractional derivative \cite{hilfer00,kilbas2006theory}, the latter is usually the most studied since the Caputo fractional derivative allows us a physical interpretation to problems with initial conditions, this derivative fulfills the property of the classical calculus that the derivative of a constant is null regardless of the order $\alpha$ of the derivative, however this does not occur with the R-L fractional derivative, and this characteristic can be used to solve nonlinear systems \cite{torres2020reduction,torres2020fractional,torres2020nonlinear}.

Unlike the Caputo fractional derivative, the R-L fractional derivative does not allow for a physical interpretation to the problems with initial condition because its use induces  fractional initial conditions, however the fact that this derivative does not cancel the constants for $\alpha$, with $\alpha\notin \nset{N}$, allows to obtain a \com{spectrum} of the behavior of the constants for different orders of the derivative, which is not possible with conventional calculus.

\subsection{Introduction to the Riemann-Liouville Fractional Derivative}

One of the key pieces in the study of fractional calculus is the iterated integral, which is defined as follows \cite{hilfer00}:

\begin{definition}
Let $ L_{loc} ^ 1 (a, b) $, the space of locally integrable functions in the interval $ (a, b) $. If $ f $ is a function such that $ f \in L_ {loc} ^ 1 (a, \infty) $, then the $n$-th iterated integral of the function $ f $ is given by 

\begin{eqnarray}\label{eq:3-001}
\begin{array}{c}
\ds \ifr{}{a}{I}{x}{n} f(x)=\ifr{}{a}{I}{x}{}\left(\ifr{}{a}{I}{x}{n-1} f(x)  \right)=\frac{1}{(n-1)!}\int_a^x(x-t)^{n-1}f(t)dt,
\end{array}
\end{eqnarray}

where

\begin{eqnarray*}
\ifr{}{a}{I}{x}{} f(x):=\int_a^x f(t)dt.
\end{eqnarray*}

\end{definition}

Considerate that $ (n-1)! = \gam{n} $
, a generalization of \eqref{eq:3-001} may be obtained for an arbitrary order $ \alpha> 0 $

\begin{eqnarray}\label{eq:3-002}
\ifr{}{a}{I}{x}{\alpha} f(x)=\dfrac{1}{\gam{\alpha}}\int_a^x(x-t)^{\alpha-1}f(t)dt,
\end{eqnarray}

similarly, if $ f \in L_{loc} ^ 1 (- \infty, b) $, we may define

\begin{eqnarray}\label{eq:3-003}
\ifr{}{x}{I}{b}{\alpha} f(x)=\dfrac{1}{\gam{\alpha}}\int_x^b(t-x)^{\alpha -1}f(t)dt,
\end{eqnarray} 

the equations \eqref{eq:3-002} and \eqref{eq:3-003}, correspond to the definitions of \textbf{right and left Riemann-Liouville fractional integral}, respectively. The fractional integrals fulfill the  \textbf{semigroup property}, which is given in the following proposition \cite{hilfer00}:

\begin{proposition}
Let $ f $ be a function. If $ f \in L_{loc} ^ 1 (a, \infty) $, then the fractional integrals of $ f $ fulfill that

\begin{eqnarray}\label{eq:3-004}
\ifr{}{a}{I}{x}{\alpha} \ifr{}{a}{I}{x}{\beta}f(x) = \ifr{}{a}{I}{x}{\alpha + \beta}f(x),& \alpha,\beta>0.
\end{eqnarray}

\end{proposition}

From the previous result, and considering that the operator $ d / dx $  is the inverse operator to the left of the operator $ \ifr {}{a}{I}{x}{} $, any integral $ \alpha$-th of a function $ f \in L_{loc} ^ 1 (a, \infty) $ may be written as

\begin{eqnarray}\label{eq:3-005}
\ifr{}{a}{I}{x}{\alpha}f(x)=\dfrac{d^n}{dx^n}\left( \ifr{}{a}{I}{x}{n}\ifr{}{a}{I}{x}{\alpha}f(x) \right)=\dfrac{d^n}{dx^n}\left( \ifr{}{a}{I}{x}{n+\alpha}f(x)\right).
\end{eqnarray}
 
Considering \eqref{eq:3-002} and \eqref{eq:3-005}, we can build the \textbf{(right) Riemann-Liouville fractional derivative} as follows \cite{hilfer00,kilbas2006theory}:

\begin{eqnarray}\label{eq:3-006}
\normalsize
\begin{array}{c}
\ifr{}{a}{D}{x}{\alpha}f(x) := \left\{
\begin{array}{cc}
\ds \ifr{}{a}{I}{x}{-\alpha}f(x), &\mbox{if }\alpha<0\\  
\ds \dfrac{d^n}{dx^n}\left( \ifr{}{a}{I}{x}{n-\alpha}f(x)\right), & \mbox{if }\alpha\geq 0
\end{array}
\right. ,
\end{array}
\end{eqnarray}

where  $ n = \lceil \alpha \rceil$.

\subsubsection{Construction of the Riemann-Liouville Fractional Derivative}

The  R-L fractional derivative is constructed in a simplified way, taking into account that the integral operator is defined for a locally integrable function $ f$, that is,  $f\in L_ {loc} ^ 1 (a, \infty)$, then

\begin{eqnarray*}
\ifr{}{a}{I}{x}{}f(x)=\int_a^xf(t)dt,
\end{eqnarray*}

applying two times the integral operator

\begin{eqnarray*}
\ifr{}{a}{I}{x}{2}f(x)= \int_a^x \left( \int_a^{x_1}f(t)dt\right)dx_1=\int_a^x \left(\ifr{}{a}{I}{x_1}{}f(x_1) \right)dx_1,
\end{eqnarray*}

doing an integration by parts, taking $u=\ifr{}{a}{I}{x_1}{}f(x_1)$ and $dv=dx_1$, as a consequence

\begin{align}
\ifr{}{a}{I}{x}{2}f(x)=& x_1 \ifr{}{a}{I}{x_1}{}f(x_1) \bigg{|}_{a}^x-\int_a^x x_1 f(x_1)dx_1 \nonumber \\
=&x \ifr{}{a}{I}{x}{}f(x)- \ifr{}{a}{I}{x}{}\left(x f(x)\right) \nonumber\\
=& \int_a^x (x-t)f(t)dt,
\end{align}

repeating the previous process, applying three times the
integral operator

\begin{eqnarray*}
\ifr{}{a}{I}{x}{3}f(x)= \int_a^x \left(\ifr{}{a}{I}{x_1}{2}f(x_1) \right)dx_1,
\end{eqnarray*}

doing an integration by parts, taking $u=\ifr{}{a}{I}{x_1}{2}f(x_1)$ and $dv=dx_1$, then

\begin{align}
\ifr{}{a}{I}{x}{3}f(x)=& x_1 \ifr{}{a}{I}{x_1}{2}f(x_1) \bigg{|}_{a}^x-\int_a^x \left( x_1 \ifr{}{a}{I}{x_1}{} f(x_1)\right) dx_1 \nonumber \\
=&x \ifr{}{a}{I}{x}{2}f(x)- \ifr{}{a}{I}{x}{}\left(x\ifr{}{a}{I}{x}{} f(x)\right) \nonumber\\
=& \int_a^x  (x-t)\ifr{}{a}{I}{t}{}f(t)   dt,
\end{align}

doing again an integration by parts, taking $u=\ifr{}{a}{I}{t}{}f(t)$ and $dv=(x-t)dt$, as a consequence

\begin{align}
\ifr{}{a}{I}{x}{3}f(x)=& -\dfrac{1}{2}(x-t)^2  \ifr{}{a}{I}{t}{}f(t) \bigg{|}_{a}^x+ \dfrac{1}{2}\int_a^x (x-t)^2f(t)dt \nonumber\\
=&  \dfrac{1}{2}\int_a^x (x-t)^2f(t)dt.
\end{align}

Repeating the previous process, applying $n$ times the integral operator and doing $n-1$ integrations by parts, it is possible to obtain the following expression of the $n$-th iterated integral\cite{hilfer00}

\begin{eqnarray*}
\ifr{}{a}{I}{x}{n}f(x)=\dfrac{1}{(n-1)!}\int_a^x (x-t)^{n-1}f(t)dt,
\end{eqnarray*}

to make a generalization of the previous expression, it is
enough to take into account the relationship between the
gamma function and the factorial function, $\gam{n}=(n-1)!$, and doing $n\to \alpha \in \nset{R}$,  the expression for
the (right) R-L fractional integral is obtained \cite{hilfer00}

\begin{eqnarray*}
\ifr{}{a}{I}{x}{\alpha}f(x)=\dfrac{1}{\gam{\alpha}}\int_a^x (x-t)^{\alpha-1}f(t)dt,
\end{eqnarray*}

taking into account that the differential operator $(D_x=d/dx)$ is the inverse operator to the left of the integral operator $(\ifr{}{a}{I}{x}{})$ , that is,

\begin{eqnarray*}
D_x^n\left( \ifr{}{a}{I}{x}{n}f(x)\right)=\dfrac{d^n}{dx^n}\left(\ifr{}{a}{I}{x}{n}f(x)\right)=f(x),
\end{eqnarray*}

we may consider extending the previous result analogously to the fractional calculus using the expression

\begin{eqnarray*}
\ifr{}{a}{D}{x}{\alpha}f(x):=\ifr{}{a}{I}{x}{-\alpha}f(x),
\end{eqnarray*}

unfortunately, this would cause convergence problems because the gamma function is not defined for $\alpha \in \nset{Z}_{<0}$, to solve this problem, the above expression is rewritten as

\begin{align*}
\ifr{}{a}{D}{x}{\alpha}f(x)=&  \ifr{}{a}{I}{x}{-\alpha}f(x)\\
=&\dfrac{d^n}{dx^n}\big{(} \ifr{}{a}{I}{x}{n}\big{(} \ifr{}{a}{I}{x}{-\alpha}f(x) \big{)} \big{)}\\
=&\dfrac{d^n}{dx^n}\left( \ifr{}{a}{I}{x}{n-\alpha}f(x) \right),
\end{align*}

for the above expression to make sense, it is necessary to consider $n-\alpha\geq 0$, there are infinite ways that $n$ may be taken to fulfills the above condition, but the most convenient way is to consider

\begin{eqnarray*}
n=n(\alpha),
\end{eqnarray*}

considering the above, we can define the (right) R-L fractional derivative as follows

\begin{eqnarray*}
\ifr{}{a}{D}{x}{\alpha}f(x)=\dfrac{1}{\gam{n-\alpha}}\dfrac{d^n}{dx^n}\int_a^x(x-t)^{n-\alpha-1}f(t)dt, & n=\lceil \alpha \rceil, 
\end{eqnarray*}

in such a way that the previous expression fulfills that

\begin{align*}
\lim_{\alpha \to 1}\ifr{}{a}{D}{x}{\alpha}f(x) =& \lim_{\alpha \to 1}\dfrac{d^n}{dx^n}\left(  \ifr{}{a}{I}{x}{n-\alpha}f(x) \right)\\
=& \dfrac{d}{dx}\left(  \ifr{}{a}{I}{x}{0}f(x) \right)\\
=& \dfrac{d}{dx}f(x).
\end{align*}

Finally, it is possible to unify the R-L fractional operators, fractional integral and fractional derivative, and define the \textbf{(right) Riemann-Liouville fractional derivative} as follows \cite{hilfer00,kilbas2006theory}:

\begin{eqnarray*}
\normalsize
\begin{array}{c}
\ifr{}{a}{D}{x}{\alpha}f(x) := \left\{
\begin{array}{cc}
\ds \ifr{}{a}{I}{x}{-\alpha}f(x), &\mbox{if }\alpha<0\\  
\ds \dfrac{d^n}{dx^n}\left( \ifr{}{a}{I}{x}{n-\alpha}f(x)\right), & \mbox{if }\alpha\geq 0
\end{array}
\right. ,
\end{array}
\end{eqnarray*}

where  $ n = \lceil \alpha \rceil$.

\subsubsection{Examples of the Riemann-Liouville Fractional Derivative}

Before continuing, it is necessary to define the Beta function and the incomplete Beta function \cite{arfken85}, which are defined as follows

\begin{eqnarray}
B(p,q):=\int_0^1 t^{p-1}(1-t)^{q-1}dt, & \ds B_r(p,q):=\int_0^r t^{p-1}(1-t)^{q-1}dt, 
\end{eqnarray}

where $p$ and $q$ are positive values. Considering the following proposition:

\begin{proposition}\label{prop:3-001} 
Let $f$ be a function, with 

\begin{eqnarray*}
f(x)=(x-c)^\mu , & \mu>-1, & c\in \nset{R},
\end{eqnarray*}

then for all $\alpha \in \nset{R}\setminus \nset{Z}$, the Riemann-Liouville fractional derivative of the above function may be written as

\begin{eqnarray}
\normalsize
\begin{array}{c}
\ifr{}{a}{D}{x}{\alpha}f(x) = \left\{
\begin{array}{cc}
\ds \dfrac{\gam{\mu +1}}{\gam{\mu-\alpha+1}}(x-c)^{\mu-\alpha}G_{-\alpha}\left(\dfrac{a-c}{x-c},\mu+1\right), &\mbox{if }\alpha<0 \vspace{0.1cm}\\  
\ds \sum_{k=0}^n \binom{n}{k}  \dfrac{\gam{\mu +1}}{\gam{\mu +n-\alpha-k+1}}  (x-c)^{\mu +n-\alpha-k} G_{n-\alpha}^{(n-k)}\left(\dfrac{a-c}{x-c},\mu +1 \right), & \mbox{if }\alpha\geq 0
\end{array}
\right.
\end{array}, 
\end{eqnarray}

where

\begin{eqnarray}\label{eq:3-007}
G_{\alpha}\left(\dfrac{a-c}{x-c},\mu+1 \right):= 1  - \dfrac{B_{\frac{a-c}{x-c}}(\mu+1,\alpha)}{B(\mu+1,\alpha)}.
\end{eqnarray}

\begin{proof}
The Riemann-Liouville fractional derivative of the function $f(x)$, through the equation \eqref{eq:3-006}, presents two cases:

\begin{itemize}
\item[i)] If $\alpha <0$, then :

\begin{eqnarray*}
\ifr{}{a}{D}{x}{\alpha}f(x)=\dfrac{1}{\gam{-\alpha}}\int_a^x(x-t)^{-\alpha-1}(t-c)^\mu dt,
\end{eqnarray*}

taking the change of variable $t=c+(x-c)u$ in the previous expression

\begin{align*}
\ifr{}{a}{D}{x}{\alpha}f(x)=\dfrac{(x-c)^{\mu-\alpha}}{\gam{-\alpha}}\int_{\frac{a-c}{x-c}}^1(1-u)^{-\alpha-1}u^\mu du,
\end{align*}

the above result may be rewritten in terms of  the Beta function and the incomplete Beta function as follows

\begin{align*}
\ifr{}{a}{D}{x}{\alpha}f(x)
=&\dfrac{(x-c)^{\mu-\alpha}}{\gam{-\alpha}}\left( B(\mu+1,-\alpha)  -B_{\frac{a-c}{x-c}}(\mu+1,-\alpha) \right) \\
=& B(\mu +1,-\alpha)\dfrac{(x-c)^{\mu-\alpha}}{\gam{-\alpha}}\left( 1  - \dfrac{B_{\frac{a-c}{x-c}}(\mu+1,-\alpha)}{B(\mu+1,-\alpha)} \right),
\end{align*}

and considering \eqref{eq:3-007}, we obtain that

\begin{eqnarray}\label{eq:3-008}
\ifr{}{a}{D}{x}{\alpha}(x-c)^\mu= \dfrac{\gam{\mu +1}}{\gam{\mu-\alpha+1}}(x-c)^{\mu-\alpha}G_{-\alpha}\left(\dfrac{a-c}{x-c},\mu+1\right).
\end{eqnarray}

\item[ii)] If $\alpha \geq 0$, then:

\begin{eqnarray*}
\ifr{}{a}{D}{x}{\alpha}f(x)=\dfrac{1}{\gam{n-\alpha}}\dfrac{d^n}{dx^n}\int_a^x(x-t)^{n-\alpha-1}(t-c)^\mu dt,
\end{eqnarray*}

taking the change of variable $t=c+(x-c)u$ in the previous expression

\begin{align*}
\ifr{}{a}{D}{x}{\alpha}f(x)=\dfrac{1}{\gam{n-\alpha}}\dfrac{d^n}{dx^n} \left[ (x-c)^{\mu +n-\alpha} \int_{\frac{a-c}{x-c}}^1(1-u)^{n-\alpha-1}u^\mu du \right] ,
\end{align*}

the above result may be rewritten in terms of  the Beta function and the incomplete Beta function as follows

\begin{align*}
\ifr{}{a}{D}{x}{\alpha}f(x)
=&\dfrac{1}{\gam{n-\alpha}}\dfrac{d^n}{dx^n} \left[ (x-c)^{\mu+n-\alpha}\left( B(\mu +1,n-\alpha)  -B_{\frac{a-c}{x-c}}(\mu +1,n-\alpha)\right) \right] \\
=&\dfrac{B(\mu +1,n-\alpha)}{\gam{n-\alpha}}\dfrac{d^n}{dx^n} \left[ (x-c)^{\mu +n-\alpha}\left( 1  - \dfrac{B_{\frac{a-c}{x-c}}(\mu +1,n-\alpha)}{B(\mu +1,n-\alpha)} \right) \right] ,
\end{align*}

and considering \eqref{eq:3-007}, we obtain that

\begin{align*}
\ifr{}{a}{D}{x}{\alpha}f(x)=&\dfrac{\gam{\mu +1}}{\gam{\mu +n-\alpha+1}}\dfrac{d^n}{dx^n} \left[ (x-c)^{\mu +n-\alpha}G_{n-\alpha}\left(\dfrac{a-c}{x-c},\mu +1 \right) \right] \\
=&\dfrac{\gam{\mu +1}}{\gam{\mu +n-\alpha+1}}  \sum_{k=0}^n \binom{n}{k}\left(\dfrac{d^k}{dx^k} (x-c)^{\mu +n-\alpha} \right) G_{n-\alpha}^{(n-k)}\left(\dfrac{a-c}{x-c},\mu +1 \right),
\end{align*}

taking into account that in the classical calculus

\begin{eqnarray*}
\dfrac{d^k}{dx^k}(x-c)^\mu =\dfrac{\mu !}{(\mu -k)!}(x-c)^{\mu-k}=\dfrac{\gam{\mu +1}}{\gam{\mu -k+1}}(x-c)^{\mu -k},
\end{eqnarray*}

therefore

\begin{align}\label{eq:3-009}
\ifr{}{a}{D}{x}{\alpha}(x-c)^\mu
= \sum_{k=0}^n \binom{n}{k}  \dfrac{\gam{\mu +1}}{\gam{\mu +n-\alpha-k+1}}  (x-c)^{\mu +n-\alpha-k} G_{n-\alpha}^{(n-k)}\left(\dfrac{a-c}{x-c},\mu +1 \right).
\end{align}

\end{itemize}
\end{proof}

\end{proposition}

From the previous proposition, we can note that the Riemann-Liouville fractional derivative presents an explicit dependence of the value $n=\lceil \alpha \rceil$. However, there exists a particular case in which this dependence disappears, as shown in the following proposition:

\begin{proposition}\label{prop:3-002}
Let $f$ be a function, with

\begin{eqnarray*}
f(x)=(x-a)^\mu , & \mu>-1, & a\in \nset{R},
\end{eqnarray*}

then for all $\alpha\in \nset{R}\setminus \nset{Z}$, the Riemann-Liouville fractional derivative of the above function may be written in general form as

\begin{eqnarray}\label{eq:3-010}
\ifr{}{a}{D}{x}{\alpha}(x-a)^\mu=\dfrac{\gam{\mu+1}}{\gam{\mu-\alpha+x}}(x-a)^{\mu-\alpha}.
\end{eqnarray}

\begin{proof}
To prove the validity of the previous equation for all $\alpha\in \nset{R}\setminus \nset{Z}$, it is necessary to note that from the \textbf{Proposition \ref{prop:3-001}} , the following limits may be obtained

\begin{eqnarray*}
\ifr{}{a}{D}{x}{\alpha}(x-a)^\mu=\lim_{c\to a}\ifr{}{a}{D}{x}{\alpha}(x-c)^\mu,
\end{eqnarray*}

\begin{eqnarray*}
\lim_{c \to a}G_{\alpha}\left(\dfrac{a-c}{x-c},m+1 \right)=G_\alpha(0,\mu+1)=1,
\end{eqnarray*}

then consider two cases:

\begin{itemize}
\item[i)] If  $\alpha < 0$, from the equation \eqref{eq:3-008}, we obtain that

\begin{align*}
\ifr{}{a}{D}{x}{\alpha}(x-a)^\mu=&\dfrac{\gam{\mu +1}}{\gam{\mu-\alpha+1}} \lim_{c\to a} \left( (x-c)^{\mu-\alpha}G_{-\alpha}\left(\dfrac{a-c}{x-c},\mu+1\right) \right) \\
=&\dfrac{\gam{\mu +1}}{\gam{\mu-\alpha+1}} (x-a)^{\mu-\alpha}G_{-\alpha}\left(0,\mu+1\right) \\
=&\dfrac{\gam{\mu +1}}{\gam{\mu-\alpha+1}} (x-a)^{\mu-\alpha}.
\end{align*}

\item[i)] If  $\alpha \geq 0$, from the equation \eqref{eq:3-009}, we obtain that

\begin{align*}
\ifr{}{a}{D}{x}{\alpha}(x-a)^\mu=&\sum_{k=0}^n \binom{n}{k}  \dfrac{\gam{\mu +1}}{\gam{\mu +n-\alpha-k+1}}\lim_{c\to a}\left(   (x-c)^{\mu +n-\alpha-k} G_{n-\alpha}^{(n-k)}\left(\dfrac{a-c}{x-c},\mu +1 \right) \right) \\
=&\sum_{k=0}^n \binom{n}{k}  \dfrac{\gam{\mu +1}}{\gam{\mu +n-\alpha-k+1}}  (x-a)^{\mu +n-\alpha-k} G_{n-\alpha}^{(n-k)}\left(0,\mu +1 \right)  \\
=& \binom{n}{n}  \dfrac{\gam{\mu +1}}{\gam{\mu +n-\alpha-n+1}}  (x-a)^{\mu +n-\alpha-n} G_{n-\alpha}^{(0)}\left(0,\mu +1 \right)  \\
=&\dfrac{\gam{\mu +1}}{\gam{\mu-\alpha+1}} (x-a)^{\mu-\alpha}.
\end{align*}

\end{itemize}

\end{proof}

\end{proposition}

From the previous proposition, the following corollary is obtained

\begin{corollary}
Let $f:\Omega \subset \nset{R} \to \nset{R}$ be a function, with $f \in L_{loc}^1(a,\infty)$. Assuming furthermore that $f\in C^\infty(a,\infty)$,   such that $f$ may be written in terms of its Taylor series around the point $x=a$, that is,

\begin{eqnarray*}
f(x)=\sum_{k=0}^\infty \dfrac{f^{(k)}(a)}{k!}(x-a)^k,
\end{eqnarray*}

then for all $\alpha \in \nset{R}\setminus \nset{Z}$, the Riemann-Liouville fractional derivative of the aforementioned function, may be written as follows

\begin{eqnarray}
\ifr{}{a}{D}{x}{\alpha}f(x)=\sum_{k=0}^\infty \dfrac{f^{(k)}(a)}{\gam{k-\alpha+1}} (x-a)^{k-\alpha}.
\end{eqnarray}

\end{corollary}

Finally, applying the  operator \eqref{eq:3-006} with $a=0$ to the  function $ x^{\mu} $, with $\mu> -1$, from the \textbf{Proposition \ref{prop:3-002}} we obtain the following result

\begin{eqnarray}\label{eq:3-011}
\ifr{}{0}{D}{x}{\alpha}x^\mu = 
 \dfrac{\gam{\mu+1}}{\gam{\mu-\alpha+1}}x^{\mu-\alpha}, & \alpha\in \nset{R}\setminus \nset{Z}.
\end{eqnarray}

\subsection{Introduction to the Caputo Fractional Derivative}

Michele Caputo (1969) published a book and introduced a new definition of fractional derivative, he created this definition with the objective of modeling anomalous diffusion phenomena. The definition of Caputo had already been discovered independently by Gerasimov (1948). This fractional derivative is of the utmost importance since it allows us to give a physical interpretation of the initial value problems, moreover to being used to model fractional time. In some texts, it is known as the fractional derivative of Gerasimov-Caputo.

Let $ f $ be a function, such that $ f $ is $ n$-times differentiable with $ f ^{(n)} \in L_{loc}^ 1 (a, b) $, then the \textbf{(right) fractional derivative  of Caputo} is defined as \cite{kilbas2006theory}

\begin{align}\label{eq:3-012}
\ifr{C}{a}{D}{x}{\alpha}f(x):= &\ifr{}{a}{I}{x}{n-\alpha}\left( \der{d}{x}{n} f(x)\right) = \dfrac{1}{\gam{n-\alpha}}\int_{a}^{x} (x-t)^{n-\alpha -1} f^{(n)}(t)dt , 
\end{align}

where $n=\lceil \alpha \rceil$. It should be mentioned that the Caputo fractional derivative behaves as the inverse operator to the left of the Riemann-Liouville fractional integral, that is,

\begin{eqnarray*}
\ifr{C}{a}{D}{x}{\alpha}(\ifr{}{a}{I}{x}{\alpha}f(x))=f(x).
\end{eqnarray*}

On the other hand, the relation between the fractional derivatives of Caputo and Riemann-Liouville is given by the following expression \cite{kilbas2006theory}

\begin{eqnarray*}
\ifr{C}{a}{D}{x}{\alpha}f(x)=\ifr{}{a}{D}{x}{\alpha}\left(f(x)-\sum_{k=0}^{n-1}\dfrac{f^{(k)}(a)}{k!}(x-a)^k\right), 
\end{eqnarray*}

then, if $f^{(k)}(a)=0 \ \ \forall k<n$, we obtain

\begin{eqnarray*}
\ifr{C}{a}{D}{x}{\alpha}f(x)=\ifr{}{a}{D}{x}{\alpha}f(x), 
\end{eqnarray*}

considering the previous particular case, it is possible to unify the definitions of R-L fractional integral and Caputo fractional derivative as follows

\begin{eqnarray}\label{eq:3-013}
\begin{array}{c}
\ifr{C}{a}{D}{x}{\alpha}f(x) := \left\{
\begin{array}{cc}
\ds \ifr{}{a}{I}{x}{-\alpha}f(x), &\mbox{if }\alpha<0\\  
\ds \ifr{}{a}{I}{x}{n-\alpha}\left( \der{d}{x}{n} f(x)\right) , & \mbox{if }\alpha\geq 0
\end{array}
\right. .
\end{array}
\end{eqnarray}

\section{Fractional Newton-Raphson Method}

Let $\nset{P}_n(\nset{R})$ be the space of polynomials of degree less than or equal to $ n\in \nset{N} $ with real coefficients. The N-R method is characterized by the fact that if it generates divergent sequences of complex numbers, they may lead to the creation of a fractal \cite{tatham}. On the other hand, the order of the fractional derivatives seems to be closely related to the fractal dimension \cite{brambila2017fractal}, based on the above, a method was developed that makes use of the N-R method and the fractional derivatives. The N-R method is useful for finding the roots of a function $ f\in \nset{P}_n(\nset{R})$. However, this method is limited because it cannot find roots $ \xi \in \nset{C} \setminus \nset {R} $, if the sequence $ \set{x_i}_{i = 0} ^ \infty $ generated by \eqref{eq:2-005} has an initial condition $ x_0 \in \nset{R} $. To solve this problem and develop a method that has the ability to find roots, both real and complex, of a polynomial if the initial condition $ x_0 $ is real, we propose a new method, which consists of the Newton-Raphson method with the implementation of the fractional derivatives. Before continuing, it is necessary to define the \textbf{fractional Jacobian matrix} of a function $ f: \Omega \subset \nset{R} ^ n \to \nset {R}^ n $ as follows

\begin{eqnarray}\label{eq:4-001}
f^{(\alpha)}(x):=\left([f]^{(\alpha)}_{jk}(x) \right) ,
\end{eqnarray}

where

\begin{eqnarray*}
[f]^{(\alpha)}_{jk}= \partial_k^\alpha [f]_j(x):= \der{\partial}{[x]_k}{\alpha}[f]_j(x), &1\leq j,k\leq n.
\end{eqnarray*}

with $[f]_j:\nset{R}^n \to \nset{R}$. The  operator $ \partial^ \alpha / \partial [x]_k^\alpha $ denotes any fractional derivative, applied only to the variable $ [x]_k $, which fulfills the following condition of continuity respect to the order of the derivative

\begin{eqnarray*}
\lim_{\alpha \to 1}\der{\partial}{[x]_k}{\alpha}[f]_j(x)=\dfrac{\partial}{\partial [x]_k}[f]_j(x), & 1\leq j,k\leq n,
\end{eqnarray*}

then, the matrix \eqref{eq:4-001} fulfills that

\begin{eqnarray}\label{eq:4-002}
\lim_{\alpha\to 1}f^{(\alpha)}(x)=f^{(1)}(x),
\end{eqnarray}

where $ f ^{(1)} (x) $ denotes the Jacobian matrix of the  function $ f $. Considering a function $\Phi:(\nset{R}\setminus \nset{Z})\times \nset{C}^n \to \nset{C}^n$, then using as a basis the idea of the N-R method \eqref{eq:2-005}, and considering any fractional derivative that fulfills the condition \eqref{eq:4-002}, we can define the \textbf{Fractional Newton-Raphson Method} as follows \cite{torres2021fractional,torreshern2020}:

\begin{eqnarray}\label{eq:4-003}
x_{i+1}:= \Phi\left( \alpha, x_i \right)=x_i-\left(f^{\left(\alpha\right)}(x_i) \right)^{-1}  f(x_i),& i=0,1,2,\cdots.
\end{eqnarray}

For the above expression to make sense, due to the part of the integral operator that fractional derivatives usually have, and  that the F N-R method can be used in a wide variety of functions \cite{torreshern2020}, we consider in the expression \eqref{eq:4-003} that each fractional derivative is obtained for a real variable $[x]_k$, and if the result allows it, this variable is subsequently substituted by a complex variable $[x_i]_k$, that is,

\begin{eqnarray}
f^{\left(\alpha\right)}(x_i):= f^{(\alpha)}(x)\bigg{|}_{ x\longrightarrow x_i}  , & x\in \nset{R}^n, & x_i\in \nset{C}^n. 
\end{eqnarray}

\subsection{Convergence of the Fractional Newton-Raphson Method}

It should be mentioned that in general, in the F N-R method  $\norm{\Phi^{(1)}(\alpha,\xi)}\neq 0$ if $\norm{f(\xi)}=0$, and from \textbf{Corollary \ref{cor:1-001}},  \textbf{Proposition \ref{prop:2-001}}, \textbf{Proposition \ref{prop:2-002}} and the condition \eqref{eq:4-002}, any sequence $ \set{x_i} _ {i = 0} ^ \infty $ generated by the iterative method \eqref{eq:4-003} has an order of convergence at least linear, that is, the F N-R method, considering the \textbf{Theorem \ref{teo:1-001}}, may fulfill an equation analogous to the equation \eqref{eq:2-006} with $p\geq 1$,  which becomes more evident when considering $\alpha\in [1-\epsilon,1+\epsilon]\setminus \set{1}$.  The aforementioned, for the case in one dimension, may be observed in the following proposition:

\begin{proposition}

Let $f: \Omega \subset \nset{R} \to \nset{R}$ be a function with a zero $\xi \in \Omega$. Then any sequence $\set{x_i}_{i=0}^\infty$ generated by the iteration function of the F N-R method, such that $x_i \to \xi$, fulfills the following condition:

\begin{eqnarray}\label{eq:4-004}
\abs{x_{i+1}-\xi}\leq \dfrac{\abs{\Phi^{(p)}(\alpha,\xi)}}{p!} \abs{x_i-\xi}^p,
\end{eqnarray}

where

\begin{eqnarray}
p=\left\{\begin{matrix}
1, & \mbox{ if } \alpha \neq 1  \mbox{ and }f^{(\alpha)}(\xi)\neq 0 \vspace{0.1cm}\\
2, & \mbox{ if } \alpha = 1 \mbox{ and }f^{(1)}(\xi)\neq 0
\end{matrix}\right. .
\end{eqnarray}

\begin{proof}\label{prop:4-001}

Considering the iteration function of the F N-R method

\begin{eqnarray*}
\Phi(\alpha,x)=x-\left(f^{(\alpha)}(x) \right)^{-1}f(x),
\end{eqnarray*}

and calculating its first and second derivative

\begin{eqnarray*}
\Phi^{(1)}(\alpha,x)=1-\left(f^{(\alpha)}(x) \right)^{-1}f^{(1)}(x) +f(x) \left[ \left(f^{(\alpha)}(x) \right)^{-2} D_x f^{(\alpha)}(x)\right]  ,
\end{eqnarray*}

\begin{align*}
\Phi^{(2)}(\alpha,x)=&
f(x)\left[ \left(f^{(\alpha)}(x) \right)^{-2}D_x^2 f^{(\alpha)}(x) - 2  \left(f^{(\alpha)} (x)\right)^{-3}\left(D_x f^{(\alpha)}(x) \right)^2 \right]\\
&+ 2 \left(f^{(\alpha)}(x) \right)^{-2} f^{(1)}(x) D_xf^{(\alpha)}(x)- \left(f^{(\alpha)}(x) \right)^{-1} f^{(2)}(x),
\end{align*}

then, assuming that $f^{(\alpha)}(\xi)\neq 0 \ \forall \alpha \in (\nset{R}\setminus \nset{Z}) \cup \set{1}$, and taking into account the condition  \eqref{eq:4-002} together with the fact that $\xi$ is a zero of $f$, we obtain that

\begin{eqnarray*}
\lim_{x \to \xi} \Phi(\alpha,x)=\xi,
\end{eqnarray*}

\begin{eqnarray*}
\lim_{x\to \xi} \Phi^{(1)}(\alpha, x)= \left\{
\begin{array}{cc}
1-\left(f^{(\alpha)}(\xi) \right)^{-1}f^{(1)}(\xi) &, \mbox{ if } \alpha \neq 1 \vspace{0.1cm}\\
0 &, \mbox{ if } \alpha = 1
\end{array}\right.,
\end{eqnarray*}

\begin{eqnarray*}
\lim_{x\to \xi} \Phi^{(2)}(\alpha, x)= \left\{
\begin{array}{cc}
2 \left(f^{(\alpha)}(\xi) \right)^{-2} f^{(1)}(\xi) D_xf^{(\alpha)}(\xi)- \left(f^{(\alpha)}(\xi) \right)^{-1} f^{(2)}(\xi) &, \mbox{ if } \alpha \neq 1 \vspace{0.1cm}\\
\left(f^{(\alpha)}(\xi) \right)^{-1} f^{(2)}(\xi) &, \mbox{ if } \alpha = 1
\end{array}\right.,
\end{eqnarray*}

as a consequence, from the \textbf{Theorem \ref{teo:1-001}}, the F N-R method has an order of convergence at least linear, that is, fulfills the equation \eqref{eq:4-004} with $p\geq 1$.

\end{proof}

\end{proposition}

From the above proposition, together with the \textbf{Proposition \ref{prop:2-001}}, it may be obtained that almost any fractional iterative method that has a similar structure to the fractional Newton-Raphson method\cite{akgul2019fractional,cordero2019variant,gdawiec2019visual,gdawiec2020newton,torreshern2020}, has the ability to change from an order of convergence (at least) linear to an order of convergence (at least) quadratic,  as long as the method fulfills the condition \eqref{eq:4-002}.   An alternative to achieve the change in the order of convergence of some fractional iterative method, analogous to F N-R method, is to replace the constant value $\alpha$ in the order of the fractional derivatives by some function that guarantees that the condition \eqref{eq:4-002} is fulfilled, that is,

\begin{eqnarray}\label{eq:4-005}
\alpha \in \nset{R}\setminus \nset{Z} & \longrightarrow & \alpha(x): \nset{C}^n \to (\nset{R}\setminus \nset{Z})\cup \set{1}.
\end{eqnarray}

It is necessary to mention that an example of the aforementioned may be found in the \textbf{Fractional Newton  Method},   which is defined as follows \cite{torreshern2020}:

\begin{eqnarray}\label{eq:4-006}
x_{i+1}:=\Phi(\alpha,x_i)=x_i-\left(  N_{\alpha_f }(x_i)  \right)^{-1}f(x_i), & i=0,1,2,\cdots,
\end{eqnarray}

where $ N_{\alpha_f} (x_i) $ is given by the following matrix 

\begin{eqnarray}
N_{\alpha_f}(x_i):=\left([N_{\alpha_f}]_{jk}(x_i) \right)=\left(\partial_k^{\alpha_f([x_i]_k,x_i)} [f]_j(x_i) \right).
\end{eqnarray}

with $\delta>0$, and $\alpha_{f}([x_i]_k, x_i)$ a function defined as follows

\begin{eqnarray}\label{eq:4-007}
\alpha_{f}([x_i]_k, x_i):=\left\{
\begin{array}{cc}
\alpha, & \mbox{ if } \sqrt{[x_i]_k [\overline{ x_i}]_k}\neq 0   \mbox{ and } \norm{f(x_i)}\geq \delta  \vspace{0.1cm}\\
1, & \mbox{ if } \sqrt{[x_i]_k [\overline{ x_i}]_k}=0 \mbox{ or }\norm{f(x_i)}<\delta
\end{array}\right. ,
\end{eqnarray}

the difference between the methods \eqref{eq:4-003} and \eqref{eq:4-006},  is that just for the second method may there exist a value $\delta> 0 $, such that if the sequence $ \set{x_i}_{i = 0} ^ \infty $ generated by \eqref{eq:4-006} converges to a zero $ \xi $ of $ f $, there exists a value $ k> 0 $ such that $ \forall i \geq k $, from  \textbf{Proposition \ref{prop:2-001}}, \textbf{Proposition \ref{prop:4-001}} and  condition \eqref{eq:4-002}, the sequence  has an order of convergence (at least) quadratic in $B(\xi;\delta)$.

\section{The Aitken's Method}

Due not all fractional iterative methods fulfill the condition \eqref{eq:4-002}, since not all methods have a similar structure to the F N-R method \cite{torreshern2020,torres2020reduction,torres2020fractional} , an alternative such as that of equation \eqref{eq:4-005} to accelerate its order of convergence would not be suitable. However, an alternative that may be used in general in any fractional iterative method to accelerate its convergence, is to combine the method with the \textbf{Aitken's method} \cite{stoer2013,nievergelt1991aitken}.

The Aitken's method or also known as the $ \Delta ^ 2 -$ method of Aitken \cite{stoer2013}, is one of the first and simplest methods to accelerate the convergence of a given convergent sequence $\set{x_i}_{i=0}^\infty$, that is,

\begin{eqnarray*}
\lim_{i \to \infty}\norm{x_i -\xi} \to 0,
\end{eqnarray*}

this method allows transforming the sequence $ \set{x_i}_ {i = }^\infty $ to a sequence $ \set{x_i '} _ {i = 0} ^ {\infty} $, which generally converges faster point $ \xi $ that the original sequence, Under certain circumstances, the Aitken's method may accelerate the convergence of a method that has an order of convergence  (at least) linear to an order of convergence almost quadratic, then it is generally used to accelerate the iterative methods used to find the zeros of a function  $f$ \cite{plato2003concise,stoer2013,burden2002analisis}.

To illustrate the Aitken's method for the case in one dimension, suppose that the sequence
$ \set{x_i} _ {i = 0} ^ \infty $ converges to the point $ \xi $ as a geometric sequence with factor $ k $ such that $\abs{k}<1$, that is,

\begin{eqnarray}\label{eq:5-001}
x_{i+1}- \xi = k \left(x_i-\xi \right), & i=0,1, 2,\cdots,
\end{eqnarray}

where the value of $ \xi $ may be determined using the following system of equations

\begin{align}
x_{i+1}- \xi &= k\left(x_{i}- \xi \right), \label{eq:5-002}\\
x_{i+2}- \xi &= k\left(x_{i+1}- \xi \right), \label{eq:5-003}
\end{align}

subtracting the equation \eqref{eq:5-002} from the equation \eqref{eq:5-003} we obtain the value of $k$

\begin{eqnarray*}
k= \dfrac{x_{i+2}-x_{i+1}}{x_{i+1}-x_i },
\end{eqnarray*}

placing $ \xi $ on the left side of the equation \eqref{eq:5-002}

\begin{align*}
\xi=& \dfrac{ k x_i -x_{i+1}}{k-1}\\
=&
\dfrac{ (k-1+1) x_i -x_{i+1}}{k-1} \\
=&x_i- \dfrac{ x_{i+1} -x_{i}}{k-1},
\end{align*}

and substituting the value of $k$ in the previous expression

\begin{align*}
\xi
=&x_i- \dfrac{\left( x_{i+1} -x_{i} \right) \left(x_{i+1}-x_i \right)}{\left(x_{i+2}-x_{i+1} \right) -\left(x_{i+1}-x_i \right)}\\
=&x_i- \dfrac{\left( x_{i+1} -x_{i} \right)^2}{x_{i+2}-2x_{i+1} +x_i},
\end{align*}

defining the difference operator

\begin{eqnarray*}
\Delta x_i := x_{i+1}-x_i,
\end{eqnarray*}

then

\begin{align*}
\Delta (\Delta x_i) =& \Delta^2 x_i\\
 =& \Delta x_{i+1}- \Delta x_i \\
=& x_{i+2}-2x_{i+1}+x_i,
\end{align*}

therefore, we obtain that the value of $\xi$ is given by the following expression

\begin{eqnarray}\label{eq:5-004}
\xi = x_i - \dfrac{\left( \Delta x_i\right)^2}{\Delta^2 x_i},
\end{eqnarray}

the Aitken´s method is considered taking into account the equation \eqref{eq:5-004}. The $ \Delta^2-$ method of Aitken consists in generating a new sequence $\set{x_i'}_{i=0}^\infty$, where

\begin{eqnarray}\label{eq:5-005}
x_i'=x_i- \dfrac{\left( x_{i+1} -x_{i} \right)^2}{x_{i+2}-2x_{i+1} +x_i}, 
\end{eqnarray}

such that

\begin{eqnarray*}
\lim_{i\to \infty}\abs{x_i'-\xi }\to 0.
\end{eqnarray*}

On the other hand, to note that the sequence $ \set{x_i '} _ {i = 0} ^ {\infty} $ converges more quickly to value $\xi $ than the sequence $ \set{x_i} _ {i = 0 } ^\infty $, consider the following proposition:

\begin{proposition}\label{prop:5-001}
Let $\set{x_i}_{i=0}^\infty$ be a sequence, such that $x_i \to \xi$. Then, the sequence $\set{x_i'}_{i=0}^\infty$ generated by the Aitken's method, given by \eqref{eq:5-005},  has a speed of convergence greater than the original sequence.

\begin{proof}

Suppose for the equation \eqref{eq:5-001} that the $k$ value fulfills the following conditions

\begin{eqnarray*}
\begin{array}{ccc}
k= k_0+ \delta_i  , & \displaystyle{\lim_{i \to \infty}}\delta_i= 0, & \abs{k}<1,
\end{array}
\end{eqnarray*}

then from equation \eqref{eq:5-001}

\begin{align}\label{eq:5-006}
x_{i+1}-x_i =& (x_{i+1}-\xi)-(x_i-\xi) \nonumber  \\
=& (k-1)(x_i-\xi),
\end{align}

analogously

\begin{align*}
x_{i+2}-x_{i+1} 
=& (k-1)(x_{i+1}-\xi) \\
=& (k-1)(x_{i+1}- x_i )+(k+1)(x_i-\xi) \\
=& \left[ (k-1)^2+(k+1) \right](x_i-\xi), 
\end{align*}

whereby

\begin{align}\label{eq:5-007}
x_{i+2}-2x_{i+1}+x_i 
=& (k-1)^2(x_i-\xi) \nonumber\\
=& \left[(k_0-1)^2 + \mu_i \right] (x_i-\xi) ,
\end{align}

where

\begin{eqnarray*}
\lim_{i \to \infty} \mu_i = 0,
\end{eqnarray*}

finally substituting the equations \eqref{eq:5-006} and \eqref{eq:5-007} in the equation \eqref{eq:5-005}, we obtain that 

\begin{eqnarray*}
x_i'-\xi = (x_i-\xi) - \dfrac{\left[ (k_0-1 + \delta_i)(x_i-\xi)\right]^2 }{\left[ (k_0-1)^2+ \mu_i  \right](x_i-\xi)},
\end{eqnarray*}

then

\begin{eqnarray*}
\dfrac{x_i'-\xi}{x_i-\xi}= 1- \dfrac{(k_0-1+\delta_i)^2}{(k_0-1)^2+ \mu_i},
\end{eqnarray*}

therefore

\begin{eqnarray*}
\lim_{i\to \infty}\dfrac{x_i'-\xi}{x_i-\xi} =0,
\end{eqnarray*}

which shows that in general, the speed of convergence of the sequence $ \set{x_i '} _ {i = 0} ^ {\infty} $ is greater than that of the original sequence.

\end{proof}

\end{proposition}

From the above proposition it follows that any fractional iterative method, given by the following expression

\begin{eqnarray}
x_{i+1}:=\Phi(\alpha,x_i), & i=0,1,2,\cdots,
\end{eqnarray}

may accelerate its speed of convergence using the Aitken's method, giving rise to the \textbf{Fractional Steffensen's Method}, which is defined as follows

\begin{eqnarray}\label{eq:5-008}
x_{i+1}:=\Psi(\alpha,x_i), & i=0,1,2,\cdots,
\end{eqnarray}

where 

\begin{eqnarray}
\Psi(\alpha,x_i):=x_i- \dfrac{\left( \Phi(\alpha,x_i) -x_{i} \right)^2}{\Phi\left(\alpha, \Phi(\alpha,x_i) \right)-2\Phi(\alpha,x_i) +x_i}.
\end{eqnarray}

It is necessary to mention that the fractional iterative method \eqref{eq:5-008} may be extended to the multimensional case, since the Aitken's method is also defined for the case of several variables \cite{nievergelt1991aitken}.

\subsection{Results of the Fractional Newton-Raphson Method with the Aitken's Method}

Examples of the implementation of the F N-R method and the Aitken's method for the multidimensional case may be found in the references \cite{torreshern2020} and \cite{nievergelt1991aitken}, respectively. However, to maintain an illustrative character, the following examples are solved for the case in one dimension and using the R-L fractional derivative \eqref{eq:3-006}. Instructions for implementing the F N-R method, along with information to provide values $\alpha \in [0.7,1.3]\setminus \set{1}$ are found in the reference \cite{torreshern2020}.  For rounding reasons, for the examples the following function is defined

\begin{eqnarray}\label{eq:5-009}
\rnd{x}{m}:=\left\{
\begin{array}{cc}
\re{x},& \mbox{ if \hspace{0.1cm}} \abs{\im{x}}\leq 10^{-m}\vspace{0.1cm}\\
x,& \mbox{ if \hspace{0.1cm}} \abs{\im{x}}> 10^{-m}\vspace{0.1cm}\\
\end{array}\right..
\end{eqnarray}

Combining the function \eqref{eq:5-009} with the methods \eqref{eq:4-003} and \eqref{eq:5-008}, the following iterative methods are defined

\begin{align}
x_{i+1}:=\rnd{\Phi(\alpha, x_i)}{5}, &\hspace{0.2cm} i=0,1,2\cdots, \vspace{0.1cm} \label{eq:5-010}\\
x_{i+1}:=\rnd{\Psi(\alpha, x_i)}{5}, &\hspace{0.2cm} i=0,1,2\cdots. \label{eq:5-011}
\end{align}

\begin{example}

Let $f$ be a function, with

\begin{align*}
f(x)=&-64.23x^{14}-72.74x^{13}-61.66x^{12}+32.26x^{11}\\
&+32.3x^{10}-41.37x^9+20.18x^8+4.32x^7\\
&-5.67x^6+17.41x^5-78.6x^4-48.27x^3\\
&
-19.31x^2+77.92x-45.03.
\end{align*}

Then the initial condition $x_0=0.68$ is chosen to use the iterative methods  given by \eqref{eq:5-010} and \eqref{eq:5-011}. Consequently, we obtain the results of Table \ref{tab:01} and Table \ref{tab:011}.

\begin{itemize}

\item[•] F N-R method without Aitken's method

\begin{table}[!ht]
\centering
\footnotesize
$
\begin{array}{c|ccccc}
\toprule
&\alpha& x_n&\norm{x_n - x_{n-1} }_2  &\norm{f\left(x_n \right)}_2& n \\ 
\midrule
    1     & 0.91615 & 0.89785325 - 0.29205146i & 2.51794E-07 & 4.71496E-05 & 33 \\
    2     & 0.92801 & 0.51558282 - 0.33422358i & 3.48673E-06 & 6.66996E-05 & 24 \\
    3     & 1.02258 & 0.51558313 + 0.33422363i & 3.71927E-06 & 2.40337E-05 & 39 \\
    4     & 1.05068 & 0.89785316 + 0.29205147i & 1.90000E-07 & 2.09988E-05 & 67 \\
    5     & 1.08961 & -0.67766754 - 0.6665907i & 1.94165E-07 & 4.82585E-05 & 93 \\
\bottomrule
\end{array}
$
\caption{Results obtained using the iterative method \eqref{eq:5-010}.}\label{tab:01}
\end{table}

\item[•] F N-R method with Aitken's method

\begin{table}[!ht]
\centering
\footnotesize
$
\begin{array}{c|ccccc}
\toprule
&\alpha& x_n&\norm{x_n - x_{n-1} }_2  &\norm{f\left(x_n \right)}_2& n \\ 
\midrule
    1     & 0.75183 & -0.07482891 + 1.01883531i & 5.97858E-05 & 5.16000E-05 & 7 \\
    2     & 0.76212 & 0.5349748 + 0.82703357i & 7.87711E-05 & 9.41344E-05 & 7 \\
    3     & 0.81497 & -0.6673652 - 1.16572645i & 1.30602E-05 & 3.96105E-05 & 9 \\
    4     & 0.90872 & 0.51558331 + 0.33422363i & 1.87279E-05 & 6.68077E-07 & 6 \\
    5     & 0.91165 & -0.6673652 + 1.16572645i & 2.14264E-06 & 3.96105E-05 & 7 \\
    6     & 0.91319 & -1.09479585 - 0.25179056i & 3.54154E-06 & 2.03933E-05 & 8 \\
    7     & 0.91350 & 0.53497473 - 0.82703357i & 7.19173E-06 & 3.77200E-06 & 8 \\
    8     & 0.91615 & 0.89785319 - 0.29205147i & 1.12395E-05 & 2.14471E-06 & 5 \\
    9     & 0.91722 & -0.67766754 + 0.66659064i & 2.18453E-05 & 2.08269E-06 & 6 \\
    10    & 0.91735 & -1.09479584 + 0.25179057i & 1.90520E-06 & 1.86462E-05 & 6 \\
    11    & 0.92801 & 0.89785318 + 0.29205147i & 5.94379E-05 & 6.07221E-06 & 6 \\
    12    & 1.02258 & -0.07482894 - 1.01883532i & 7.24315E-05 & 5.81437E-06 & 8 \\
    13    & 1.05068 & 0.51558331 - 0.33422363i & 3.10364E-05 & 6.68077E-07 & 6 \\
    14    & 1.08961 & -0.67766754 - 0.66659064i & 9.16106E-06 & 2.08269E-06 & 7 \\
\bottomrule
\end{array}
$
\caption{Results obtained using the iterative method \eqref{eq:5-011}.}\label{tab:011}
\end{table}
\end{itemize}
\end{example}

\begin{example}

Let $f$ be a function, with

\begin{align*}
f(x)=&-96.98x^{15}-96.82x^{14}-3.87x^{13}+25.78x^{12}\\
&+90.68x^{11}+48.05x^{10}+50.54x^9-5.16x^8\\
&+47.01x^7+90.23x^6+87.09x^5+53.09x^4\\
&+15.38x^3+97.98x^2-61.98x+14.69.
\end{align*}

Then the initial condition $x_0=0.15$ is chosen to use the iterative methods  given by \eqref{eq:5-010} and \eqref{eq:5-011}. Consequently, we obtain the results of  Table \ref{tab:02} and Table \ref{tab:022}.

\begin{itemize}
\item[•] F N-R method without Aitken's method

\begin{table}[!ht]
\centering
\footnotesize
$
\begin{array}{c|ccccc}
\toprule
&\alpha& x_n&\norm{x_n - x_{n-1} }_2  &\norm{f\left(x_n \right)}_2& n \\ 
\midrule
    1     & 0.89914 & 0.77399752 - 0.54762166i & 1.34164E-07 & 4.92991E-05 & 25 \\
    2     & 0.92053 & 0.28667106 - 0.19437664i & 2.73412E-06 & 2.74793E-05 & 14 \\
    3     & 0.94145 & 0.77399749 + 0.54762167i & 1.64012E-07 & 4.00860E-05 & 16 \\
    4     & 1.03241 & -1.16959783 + 0.06354746i & 3.91152E-07 & 5.57585E-05 & 24 \\
    5     & 1.03400 & 1.16397069 & 5.00000E-08 & 3.16258E-05 & 33 \\
    6     & 1.09493 & -0.8252629 + 0.64969532i & 6.40312E-08 & 4.07721E-05 & 57 \\
\bottomrule
\end{array}
$
\caption{Results obtained using the iterative method \eqref{eq:5-010}.}\label{tab:02}
\end{table}

\item[•] F N-R method with Aitken's method

\begin{table}[!ht]
\centering
\footnotesize
$
\begin{array}{c|ccccc}
\toprule
&\alpha& x_n&\norm{x_n - x_{n-1} }_2  &\norm{f\left(x_n \right)}_2& n \\ 
\midrule
    1     & 0.74272 & -1.16959778 - 0.06354744i & 3.42813E-05 & 1.46807E-05 & 5 \\
    2     & 0.74461 & 0.28667089 + 0.19437617i & 9.23734E-05 & 1.47225E-06 & 6 \\
    3     & 0.81520 & 0.3039298 - 0.77330878i & 1.32087E-05 & 1.23726E-06 & 7 \\
    4     & 0.89914 & 0.03271741 + 1.02608474i & 3.12570E-07 & 5.51104E-06 & 8 \\
    5     & 0.91170 & 0.03271742 - 1.02608474i & 5.58523E-05 & 1.10782E-05 & 7 \\
    6     & 0.92053 & -0.48361537 - 0.92838297i & 1.30801E-06 & 3.86093E-06 & 8 \\
    7     & 0.94145 & 0.28667089 - 0.1943762i & 3.82966E-05 & 2.79968E-07 & 6 \\
    8     & 1.02586 & 0.7739975 + 0.5476217i & 8.80368E-05 & 1.52509E-05 & 6 \\
    9     & 1.03241 & -0.82526289 + 0.64969532i & 3.99657E-05 & 2.10803E-05 & 7 \\
    10    & 1.03304 & 0.30392981 + 0.77330879i & 5.60892E-07 & 5.27621E-06 & 7 \\
    11    & 1.03400 & -0.82526288 - 0.64969531i & 7.27737E-05 & 1.74287E-05 & 7 \\
    12    & 1.03414 & 0.77399751 - 0.5476217i & 9.02351E-05 & 1.17943E-05 & 7 \\
    13    & 1.03641 & -1.16959779 + 0.06354744i & 5.65205E-05 & 2.95004E-06 & 6 \\
    14    & 1.04085 & -0.48361536 + 0.92838298i & 1.04827E-05 & 2.73789E-05 & 7 \\
    15    & 1.09493 & 1.16397069 & 1.10800E-05 & 3.16258E-05 & 3 \\
\bottomrule
\end{array}
$
\caption{Results obtained using the iterative method \eqref{eq:5-011}.}\label{tab:022}
\end{table}

\end{itemize}

\end{example}

\begin{example}

Let $f$ be a function, with

\begin{align*}
f(x)=&-57.62x^{16}-56.69x^{15}-37.39x^{14}-19.91x^{13}
+35.83^{12}\\
&-72.47^{11}+44.41x^{10}+43.53x^9+59.93x^8\\
&-42.9x^7-54.24x^6+72.12x^5-22.92x^4\\
&+56.39x^3+15.8x^2+60.05x+55.31.
\end{align*}

Then the initial condition $x_0=0.83$ is chosen to use the iterative methods  given by \eqref{eq:5-010} and \eqref{eq:5-011}. Consequently, we obtain the results of Table \ref{tab:03} and Table \ref{tab:033}.

\newpage

\begin{itemize}

\item[•] F N-R method without Aitken's method

\begin{table}[!ht]
\centering
\footnotesize
$
\begin{array}{c|ccccc}
\toprule
&\alpha& x_n&\norm{x_n - x_{n-1} }_2  &\norm{f\left(x_n \right)}_2& n \\ 
\midrule
    1     & 1.12775 & -0.62435238 & 1.25000E-06 & 6.46233E-05 & 4 \\
    2     & 1.16720 & -0.35983765 - 1.18135268i & 2.00000E-08 & 8.22921E-05 & 65 \\
    3     & 1.18548 & -1.00133943 & 7.00000E-08 & 6.28497E-05 & 9 \\
    4     & 1.20948 & 0.58999226 + 0.8669969i & 3.16228E-08 & 5.80021E-05 & 49 \\
    5     & 1.20951 & -0.35983763 + 1.18135266i & 1.41421E-08 & 8.64337E-05 & 83 \\
    6     & 1.21722 & -1.36995269 & 1.00000E-08 & 8.96035E-05 & 20 \\
\bottomrule
\end{array}
$
\caption{Results obtained using the iterative method \eqref{eq:5-010}.}\label{tab:03}
\end{table}

\item[•] F N-R method with Aitken's method

\begin{table}[!ht]
\centering
\footnotesize
$
\begin{array}{c|ccccc}
\toprule
&\alpha& x_n&\norm{x_n - x_{n-1} }_2  &\norm{f\left(x_n \right)}_2& n \\ 
\midrule
    1     & 0.81015 & 0.88121184 + 0.42696222i & 6.19967E-05 & 4.18492E-05 & 7 \\
    2     & 0.81393 & 0.58999222 - 0.8669969i & 5.86870E-05 & 6.29524E-05 & 6 \\
    3     & 0.86589 & 1.03423973 & 7.51161E-05 & 7.69988E-05 & 4 \\
    4     & 0.86937 & -0.70050492 - 0.785771i & 5.47139E-05 & 3.06730E-05 & 6 \\
    5     & 0.86963 & -0.28661369 - 0.80840641i & 7.38343E-05 & 3.63141E-06 & 5 \\
    6     & 0.87356 & 0.88121186 - 0.42696221i & 9.58783E-05 & 6.43810E-05 & 5 \\
    7     & 0.87663 & 0.58999225 + 0.86699686i & 3.71623E-05 & 1.57334E-05 & 6 \\
    8     & 0.89759 & -0.35983764 - 1.18135267i & 1.70880E-07 & 2.53547E-05 & 6 \\
    9     & 1.12775 & -0.62435276 & 6.99000E-06 & 6.01978E-07 & 2 \\
    10    & 1.16720 & 0.36452488 - 0.8328782i & 2.69320E-06 & 6.90389E-06 & 6 \\
    11    & 1.16881 & 0.36452487 + 0.8328782i & 5.13697E-05 & 1.07542E-05 & 5 \\
    12    & 1.18548 & -1.00133952 & 8.49000E-06 & 5.58349E-06 & 2 \\
    13    & 1.20862 & -0.35983764 + 1.18135267i & 3.14841E-06 & 2.53547E-05 & 7 \\
    14    & 1.20948 & -0.70050488 + 0.785771i & 1.93591E-05 & 4.44105E-05 & 5 \\
    15    & 1.20951 & -0.28661368 + 0.80840642i & 3.47917E-05 & 2.80317E-06 & 5 \\
    16    & 1.21722 & -1.36995269 & 5.70959E-05 & 8.96035E-05 & 2 \\
\bottomrule
\end{array}
$
\caption{Results obtained using the iterative method \eqref{eq:5-011}.}\label{tab:033}
\end{table}

\end{itemize}

\end{example}

\begin{example}

Let $f$ be a function, with

\begin{align*}
f(x)=& \sin(x)-\dfrac{x}{50},
\end{align*}

and assuming that

\begin{eqnarray*}
f^{(\alpha)}(x)\approx \ifr{}{0}{D}{x}{\alpha}\left(\sum_{k=0}^{50}\dfrac{(-1)^k}{\gam{2k+2}}x^{2k+1}-\dfrac{x}{50} \right).
\end{eqnarray*}

Then the initial condition $x_0=1.27$ is chosen to use the iterative methods  given by \eqref{eq:5-010} and \eqref{eq:5-011}. Consequently, we obtain the results of Table \ref{tab:04} and Table \ref{tab:044}.

\begin{itemize}

\item[•] F N-R method without Aitken's method

\begin{table}[!ht]
\centering
\footnotesize
$
\begin{array}{c|ccccc}
\toprule
&\alpha& x_n&\norm{x_n - x_{n-1} }_2  &\norm{f\left(x_n \right)}_2& n \\ 
\midrule
    1     & 1.09965 & -9.23893088 & 7.52000E-06 & 4.76229E-07 & 17 \\
    2     & 1.16346 & -12.8257859 & 6.95000E-06 & 2.72119E-07 & 14 \\
    3     & 1.18675 & -25.67192826 & 7.40000E-07 & 3.36984E-08 & 15 \\
    4     & 1.18853 & -27.6874723 & 3.19000E-06 & 1.00228E-06 & 14 \\
    5     & 1.19167 & -32.11337957 & 6.80000E-07 & 1.85548E-08 & 19 \\
    6     & 1.19267 & -33.81479451 & 4.92000E-06 & 1.60616E-06 & 13 \\
    7     & 1.23050 & 38.5804861 & 3.91000E-06 & 4.03819E-06 & 13 \\
\bottomrule
\end{array}
$
\caption{Results obtained using the iterative method \eqref{eq:5-010}.}\label{tab:04}
\end{table}

\newpage

\item[•] F N-R method with Aitken´s method

\begin{table}[!ht]
\centering
\footnotesize
$
\begin{array}{c|ccccc}
\toprule
&\alpha& x_n&\norm{x_n - x_{n-1} }_2  &\norm{f\left(x_n \right)}_2& n \\ 
\midrule
    1     & 0.96084 & -3.07995454 & 4.77000E-06 & 3.63773E-10 & 6 \\
    2     & 0.96089 & 3.07995454 & 4.30000E-07 & 3.63773E-10 & 6 \\
    3     & 0.96102 & 6.41177489 & 8.10000E-07 & 4.14786E-10 & 5 \\
    4     & 1.09965 & 25.6719283 & 6.25000E-06 & 1.73317E-10 & 5 \\
    5     & 1.16346 & -12.82578618 & 1.70000E-06 & 7.08794E-09 & 3 \\
    6     & 1.18675 & -25.6719283 & 1.14000E-06 & 1.73317E-10 & 8 \\
    7     & 1.18853 & -27.68747348 & 1.67000E-06 & 3.88305E-09 & 3 \\
    8     & 1.19167 & -32.1133796 & 9.68000E-06 & 3.83964E-09 & 9 \\
    9     & 1.19267 & -33.81479667 & 8.42000E-06 & 2.81606E-08 & 4 \\
    10    & 1.23050 & 38.58047959 & 1.70000E-07 & 2.74212E-08 & 5 \\
\bottomrule
\end{array}
$
\caption{Results obtained using the iterative method \eqref{eq:5-011}.}\label{tab:044}
\end{table}

\end{itemize}

\end{example}

\section{Conclusions}

In this document, in summary, the following results are presented: In \textbf{Corollary \ref{cor:1-001}}, an alternative way is obtained to demonstrate when an iterative method has an order of convergence at least linear. Considering \textbf{Proposition \ref{prop:2-001}} together with \textbf{Proposition \ref{prop:2-002}}, it is proved that Newton's method fulfills a necessary but not sufficient condition to have an order of convergence at least quadratic. In \textbf{Proposition \ref{prop:3-001}},  the radical differences that may there exist between the results of the conventional calculus and the fractional calculus when obtaining the derivative of a function are exposed, which is a consequence of dependency of the integer parameter $n(\alpha )$, which generally has the fractional derivative. In \textbf{Proposition \ref{prop:3-002}}, it is proved that under certain conditions, the results when calculating the derivative of a function in the fractional calculus are analogous to those obtained in the conventional calculus. In \textbf{Proposition \ref{prop:4-001}}, it is proved that the F N-R method has an order of convergence at least linear, but it follows that it has the ability to gradually change to an order of convergence at least quadratic as the value $\alpha$ approaches the value of one. It also follows that the change in the order of convergence in the F N-R method may be achieved by implementing a function in the order of the fractional derivatives. In \textbf{Proposition \ref{prop:5-001}}, it is proved that any succession may accelerate its speed of convergence through the implementation of Aitken's method, with which it follows that it is an ideal alternative to accelerate the speed of convergence of any fractional iterative method that does not have a structure similar to the F N-R method.

Taking into account the results in this document, although there are surely different alternatives to accelerate the  speed of convergence of the fractional iterative methods, take for example the strategy of changing the constant order $\alpha$ of the fractional derivative by a function and giving rise to the method \eqref{eq:4-006}, the Aitken's method is a simple and efficient method to accelerate the speed of convergence of any fractional iterative method, in particular for the F N-R method, due it presents an order of convergence at least linear for the case in which the order of the derivative is different from one. Then in conjunction with the Aitken method, it is concluded that the F N-R method becomes an efficient iterative method to calculate the largest possible number of zeros of a function.

\bibliography{Biblio}

\bibliographystyle{unsrt}

\end{document}